\documentclass[12pt]{article}
\usepackage{amsfonts,mathrsfs,amssymb,amsthm,mathptm}
\usepackage{amsmath,amscd}
\usepackage{mathptm,pslatex}
\usepackage{multirow}
\usepackage{color}
\usepackage[dvips]{graphicx}
\usepackage[all]{xy}
\usepackage{graphicx}
\usepackage{fancyhdr}
\usepackage{url}
\usepackage{exscale}
\usepackage{relsize}
\usepackage{bbm}

\usepackage{indentfirst}
\numberwithin{equation}{section}

\DeclareFontFamily{U}{cal}{}
\DeclareFontShape{U}{cal}{m}{n}{<->cmsy10}{}

\DeclareSymbolFont{rcal}{U}{cal}{m}{n}
\DeclareSymbolFontAlphabet{\mathcal}{rcal}

\newtheorem{Def}{Definition}[section]
\newtheorem{Bsp}[Def]{Example}
\newtheorem{Prop}[Def]{Proposition}
\newtheorem{Theo}[Def]{Theorem}
\newtheorem{Lem}[Def]{Lemma}

\theoremstyle{definition}

\newtheorem{Con}[Def]{Consruction}

\newcommand{\add}{{\rm add}}

\newcommand{\del}{{\rm del}}

\newcommand{\End}{{\rm End}}

\newcommand{\rad}{{\rm rad}}
\newcommand{\soc}{{\rm soc}}
\renewcommand{\top}{{\rm top}}
\newcommand{\pd}{{\rm pd}}

\newcommand{\fd} {{\rm fin.dim }}

\newcommand{\modcat}[1]{{\text{mod-}{#1}}}

\newcommand{\pmodcat}[1]{{\text{proj-}{#1}}}

\newcommand{\opp}{^{\rm op}}

\newcommand{\Hom}{{\rm Hom}}

\newcommand{\lra}{\longrightarrow}
\newcommand{\lraf}[1]{\stackrel{#1}{\lra}}
\newcommand{\ra}{\rightarrow}

\title{ \bf  On the asymmetry of finite delooping levels
\\
\footnotetext{Email addresses: syl13536@126.com, zhangjb@ahu.edu.cn.}
}

\author {YongLiang Sun, Jinbi Zhang
\thanks{Corresponding author} \\
{\it \scriptsize  School of Mathematics and Physics, Yancheng Institute of Technology, 224003 Yancheng, Jiangsu, P. R. China}\\
{\it \scriptsize  School of Mathematical Sciences, Anhui University, Hefei, 230601, Anhui Province, PR China}
}
\date{}
\begin{document}

\maketitle

\begin{abstract}
For any Artin algebra, we construct a related algebra that increases the delooping level on one side while decreasing it to zero on the opposite side. This dual construction corresponds to Cummings' original work on finite dimensional algebras, later extended to rings by Henning Krause. As an application, we show that the finite delooping level is not left-right symmetric.
\end{abstract}

\section{Introduction}
Let $A$ be an Artin algebra. We denote by $\mbox{mod-}A$ the category of all finitely generated right $A$-modules, and by $\pmodcat{A}$ the category of all finitely generated projective right $A$-modules. For an $A$-module $X$, we denote by $\pd(X)$ the projective dimension of $X$. The \emph{finitistic dimension} of $A$ is defined as
$$\fd(A):=\sup\{\pd(X)\mid \pd(X)<\infty,\, X\in \modcat{A}\}.$$
The finitistic dimension has attracted much interest due to the finitistic dimension conjecture, an unsolved problem asserting that the finitistic dimension of any Artin algebra is finite (see \cite{bass} or \cite[Conjecture (11), p.410]{Aus1997}).
Finitistic dimensions can be defined in terms of either left or right modules. However, the resulting value depends on this choice. In fact, the left and right finitistic dimensions of an algebra can be arbitrarily different  \cite[Example 1.2]{gkk91}.  Note that an algebra satisfies the Nakayama conjecture if and only if its opposite algebra does. 
In 1991, Happel \cite{h91} asked whether the finitistic dimension conjecture is also symmetric, that is, if an algebra has finite finitistic dimension, does its opposite algebra necessarily have finite finitistic dimension? This question has remained open for the past three decades \cite[p.27]{ch10}.
Recently, Cummings \cite{c22} showed that the finitistic dimension conjecture holds true for all finite dimensional algebras if and only if, for all finite dimensional algebras, the finitistic dimension of an algebra being finite implies that the finitistic dimension of its opposite algebra is also finite, which is true also for Artin algebras (see \cite{k22}). 
In this paper, we show a similar statement about the delooping level.
Recall that the \emph{delooping level} of a finitely generated $A$-module $X$ , as defined by Gelinas in \cite{g22}, is given by
$$\del(X):=\min\{d\mid \Omega^d(X)\mbox{ is a direct summand of } P\oplus \Omega^{d+1}(M),\mbox{ for some }P\in \pmodcat{A},\, M\in \modcat{A}\},$$
and we set $\del(X):=\infty$ if such a number $d$ does not exist, where $\Omega^d(X)$ is the $d$-th syzygy of $X$. The \emph{delooping level} of algebra $A$ is the maximum of delooping levels of simple $A$-modules:
$$\del(A):=\max\left\{\del(S)\mid S\mbox{ is simple }A\mbox{-module}\right\}.$$
Clearly, $\del(X\oplus Y)=\max\{\del(X),\del(Y)\}$ for $X,Y\in \mbox{mod-}A$, and $\del(A)=\del(\top(A_A))$, where $\top(X_A)$ denotes the top of an $A$-module $X_A$. 
It is shown in \cite{g22} that $\fd(A)\le \del(A^{\opp})$, where $A^{\opp}$ denotes the opposite algebra of $A$. Thus, the delooping level serves as a useful tool for addressing the finitistic dimension conjecture.
Recently, Xi and Zhang \cite{xi22} proved that the delooping levels are invariants under stable equivalences of algebras without nodes.
Examples of algebras with finite delooping level are abundant and capture many interesting classes of algebras, for instance, Gorenstein algebras, syzygy finite algebras, monomial algebras, radical square zero algebras and Nakayama algebras (see \cite{g22,r21,s21}).

In this paper, given any Artin algebra, we construct a related algebra and compute its delooping level. This construction is the dual of Cummings' construction, originally developed for finite dimensional algebras and later extended to all rings by Henning Krause. The main result of this paper is the following:

\begin{Theo}
\label{main-thm}
Let $A$ be an Artin algebra and $B:=T(A/\rad(A))$ the trivial extension algebra of $A/\rad(A)$ with respect to the bimodule $A/\rad(A)$. Let
$$\Lambda:=
\begin{pmatrix}
A&_AA/\rad(A)_B\\
0&B
\end{pmatrix}.$$
Then we have del$(A)\leq$ del$(\Lambda)\leq$ del$(A)+1$ and $\del(\Lambda^{\opp})=0$.
\end{Theo}

As an application, we show that the delooping level of any Artin algebra is finite if and only if, for all Artin algebras, the delooping level of an algebra being finite implies that the delooping level of its opposite algebra is also finite.
In \cite{kr23}, Kershaw and Rickard provided an example of a finite dimensional algebra with infinite delooping level. This shows that the finiteness of the delooping level is not left-right symmetric.

\section{Proof of the main result}
This section is devoted to the proof of Theorem \ref{main-thm}. Throughout the section, all algebras are Artin algebras over a commutative Artin ring $k$, and all modules are finitely generated right modules. We first recall the definition of covers of Artin algebras, as detailed in \cite{c22,k22}.

For an algebra $A$, we denote by $\rad(A)$ the Jacobson radical of $A$. For an $A$-module $X$, we denote by soc$(X)$ the largest semisimple submodule of $X$, and by top$(X)$ the largest semisimple factor module of $X$. In particular, $\top(X)\simeq X/\rad(X)$ as $A$-modules.
We consider the \emph{trivial extension}
$$T(A):=A\ltimes A^\natural$$
which is given by the bimodule $A^\natural:={_AA_A}$.
This algebra is by definition the abelian group $T(A)=A\oplus A^\natural$ with multiplication given by the formula
$$(x,y)\cdot (x',y')=(xx',xy'+yx').$$
Note that $A^\natural$ is a two-sided ideal with $T(A)/A^\natural\simeq A$ as algebras. In particular, if $A$ is a symmetric algebra, namely $A\simeq D(A)$ as $A$-$A$-bimodules, then $T(A)$ is also symmetric by \cite[Chapter IV, Proposition 3.9, p.128 ]{Aus1997}. Here $D$ stands for the usual duality of Artin algebras.

\begin{Lem}\label{T(A)}
Let $A$ be a semisimple algebra. Then
$$\rad(T(A))= A^\natural= \soc(T(A))\simeq \top(T(A)).$$
Moreover, $T(A)$ is a symmetric Nakayama algebra with $\modcat{T(A)}=\add(\top(T(A))\oplus T(A))$, and has no projective simple modules.
\end{Lem}
{\it Proof.} By \cite[Lemma 1]{k22}, we know $\rad(T(A))= A^\natural= \soc(T(A))\simeq \top(T(A)).$ Since each Artin algebra is Morita equivalent to a basic Artin algebra and their trivial extensions are Morita equivalent, we can assume that $A$ is basic. As is known, each semisimple Artin algebra is isomorphic to a product of finitely many the full matrix rings over division rings. By assumption, we write $A=\prod_{i=1}^{n}D_i$, where $D_i$ is division ring for $1\le i\le n$. Then $T(A)\simeq \prod_{i=1}^{n}T(D_i)$ as rings. Note that $T(D_i)$ is a symmetric Nakayama algebra with $\modcat{T(D_i)}=\add(\top(T(D_i))\oplus T(D_i))$, and has no projective simple modules. Thus $T(A)$ is a symmetric Nakayama algebra with $\modcat{T(A)}=\add(\top(T(A))\oplus T(A))$, and has no projective simple modules.
$\square$

\begin{Con}\label{consruction1}
(\cite{c22,k22}) Let $A$ be an Artin algebra and $B:=T(A/\rad(A))$. We view $A/\rad(A)$ as a
$B$-$A$-bimodule with left action via the canonical epimorphism
$T(A/\rad(A))\ra A/\rad(A)$,
and consider the triangular matrix algebra
$$\widetilde{A}:=
\begin{pmatrix}
A&0\\
A/\rad(A)&B
\end{pmatrix}.$$
Let $e:=\left(
\begin{smallmatrix}
1&0\\
0&0\end{smallmatrix}\right)$.
Then $\End_{\widetilde{A}}(e\widetilde{A})\simeq e \widetilde{A} e\simeq A$ as algebras. Thus $\widetilde{A}$ is said to be a \emph{cover} of $A$.
\end{Con}

\begin{Lem}\label{del=0}
$\del(\widetilde{A})=0$.
\end{Lem}
{\rm Proof.}
Notice
$\rad(\widetilde{A})
=\left(\begin{smallmatrix}
\rad(A) & 0 \\
A/\rad(A) & \rad(B)
\end{smallmatrix}\right)$.
Then $\rad^2(\widetilde{A})
=\left(\begin{smallmatrix}
\rad^2(A) & 0 \\
\rad(B)\,A/\rad(A) & \rad^2(B)
\end{smallmatrix}\right)$. By Lemma \ref{T(A)}, $\rad^2(B)=0$. It follows from the left $B$-module structure of $A/\rad(A)$ that $_BA/\rad(A)\simeq \top(_BB)$ as left $B$-modules. So, we have $\rad(B)\, A/\rad(A)=0$ and $\rad^2(\widetilde{A})
=\left(\begin{smallmatrix}
\rad^2(A) & 0 \\
0 & 0
\end{smallmatrix}\right)$.

Let $\tilde{e}:=\left(
\begin{smallmatrix}
0&0\\
0&1\end{smallmatrix}\right)$. Then $\rad^2(\tilde{e}\widetilde{A})=0$. Thus $\rad(\tilde{e}\widetilde{A})$ is semisimple and $$\rad(\tilde{e}\widetilde{A})\simeq (A/\rad(A),0,0)\oplus (0,\rad(B),0)$$ as $\widetilde{A}$-modules.
Here we identify $\widetilde{A}$-modules with triples $(X_A,Y_B,f)$, where $X$ is an $A$-module, $Y$ is a $B$-module and $f:Y\otimes_{B}A/\rad(A)\ra X$ is a homomorphism of $A$-modules. So, there exists a short exact sequence of right $\widetilde{A}$-modules
$$0
\lra (A/\rad(A),0,0)_{\widetilde{A}}\oplus (0,\rad(B),0)_{\widetilde{A}}
\lra \tilde{e}\widetilde{A}
\lraf{f} (0,\top(B),0)_{\widetilde{A}}
\lra 0.$$
Since $f$ is a projective cover of $(0,\top(B),0)_{\widetilde{A}}$, we have  $$(A/\rad(A),0,0)_{\widetilde{A}}\oplus (0,\rad(B),0)_{\widetilde{A}}\simeq \Omega((0,\top(B),0)_{\widetilde{A}}).$$
By the definition of the delooping levels, we get $\del((A/\rad(A),0,0)_{\widetilde{A}}\oplus (0,\rad(B),0)_{\widetilde{A}})=0$.
By Lemma \ref{T(A)}, $\rad(B)\simeq B/\rad(B)$ as $B$-modules. Then $\del((A/\rad(A),0,0)_{\widetilde{A}}\oplus (0,B/\rad(B),0)_{\widetilde{A}})=0$.
As is known, the simple $\widetilde{A}$-modules are either of the form $(S,0,0)$, where $S$ is a simple $A$-module, or of the form $(0,T,0)$, where $T$ is a simple $B$-module.
Thus $\del(\widetilde{A})
=\del((A/\rad(A),0,0)_{\widetilde{A}}\oplus (0,B/\rad(B),0)_{\widetilde{A}})=0$.
$\square$

Dually, we have the following construction.
\begin{Con}\label{consruction2}
Let $A$ be an Artin algebra and $B:=T(A/\rad(A))$. View $A/\rad(A)$ as an
$A$-$B$-bimodule with right action via the canonical epimorphism
$T(A/\rad(A))\ra A/\rad(A)$,
and consider the triangular matrix algebra
$$\Lambda:=
\begin{pmatrix}
A&A/\rad(A)\\
0&B
\end{pmatrix}.$$
\end{Con}
Clearly, $\Lambda^{\opp}\simeq \widetilde{A^{\opp}}$ as algebras, where $\Lambda^{\opp}$ stands for
the opposite algebra of $\Lambda$. Then we have the following direct consequence of Lemma \ref{del=0}.

\begin{Theo}\label{delopp=0}
$\del(\Lambda^{\opp})=0$.
\end{Theo}

\begin{Lem}\label{simple-seq}
Let $e':=\left(
\begin{smallmatrix}
0&0\\
0&1\end{smallmatrix}\right)\in \Lambda$. Then there exists an exact sequence
$$(\diamond)\quad 0\lra (0,B/\rad(B),0)_{\Lambda}\lra e'\Lambda \lra (0,B/\rad(B),0)_{\Lambda} \lra 0.$$ In particular, $(0,B/\rad(B),0)_{\Lambda}\simeq  \Omega((0,B/\rad(B),0)_{\Lambda})$ as $\Lambda$-modules and $\del(\top(e'\Lambda_{\Lambda}))=0$.
\end{Lem}
{\rm Proof.}  By Lemma \ref{T(A)}, there exists an exact sequence $$0\lra B/\rad(B)\lra B\lra B/\rad(B)\lra 0$$
in $\mbox{mod-}B$. It induces the following exact sequence of $\Lambda$-modules.
$$0\lra (0,B/\rad(B),0)_{\Lambda}\lra (0,B,0)_{\Lambda} \lraf{f} (0,B/\rad(B),0)_{\Lambda} \lra 0$$
It follows from $e'\Lambda=(0,B,0)_{\Lambda}$ that the desired sequence $(\diamond)$.
As $\top(e'\Lambda)\simeq (0,B/\rad(B),0)_{\Lambda}$, $f$ is a projective cover of $(0,B/\rad(B),0)_{\Lambda}$.
Thus $(0,B/\rad(B),0)_{\Lambda}\simeq  \Omega((0,B/\rad(B),0)_{\Lambda})$ as $\Lambda$-modules, and $\del(\top(e'\Lambda_{\Lambda}))=0$.
$\square$

\begin{Lem}\label{tri-del-lem}
Let $(X_A,Y_B,f)$ be a $\Lambda$-module and $M:={_AA/\rad(A)_B}$. Then there exists an exact sequence of $\Lambda$-modules
$$(\dagger)\quad 0\lra (\Omega(X_A),0,0)\oplus (0,Z_1,0)\lra (P,(P\otimes_AM)\oplus Q_B,(1_{P\otimes_AM},0)) \lraf{\gamma} (X_A,Y_B,f)_{\Lambda}\lra 0$$
such that
$\gamma$ is a projective cover of $(X_A,Y_B,f)_{\Lambda}$ with $P\in \pmodcat{A}$, $Q\in\pmodcat{B}$ and $Z_1\in\add(B/\rad(B)_B)$. Moreover, $\Omega^s((X_A,Y_B,f)_{\Lambda})\simeq (\Omega^s(X_A),0,0)_{\Lambda} \oplus (0,Z_s,0)_{\Lambda}$ with $Z_s\in\add(B/\rad(B)_B)$ for $s\ge 1$.
\end{Lem}
{\rm Proof.} Note that the simple $\Lambda$-modules are either of the form $(S,0,0)$, where $S$ is a simple $A$-module, or of the form $(0,T,0)$, where $T$ is a simple $B$-module. It follows from
$\rad(\Lambda)=
\left(\begin{smallmatrix}
\rad(A)&M\\
0&\rad(B)
\end{smallmatrix}\right)$ that
\begin{align*}
\top((X_A,Y_B,f)_{\Lambda})
&=\frac{(X_A,Y_B,f)}{\rad(\Lambda)(X_A,Y_B,f)}\\
&\simeq (\frac{X}{X\rad(A)},\frac{Y}{{\rm Im}(f)+Y\rad(B)},0)\\
&\simeq (\frac{X}{X\rad(A)},0,0)\oplus (0,\frac{Y}{{\rm Im}(f)+Y\rad(B)},0),
\end{align*}
where ${\rm Im}(f)$ stands for the image of $f$.
Let $\gamma'_1:P_A\ra X/\rad(X)$ be a projective cover of $X/\rad(X)$ and $\pi_1:X\ra X/\rad(X)$ be the natural epimorphism. Then there is a homomorphism $\gamma_1:P_A\ra X_A$ such that $\gamma'_1=\gamma_1\pi_1$ and $\gamma_1$ is a projective cover of $X_A$.
Similarly, let $\gamma'_2:Q_B\ra Y/({\rm Im}(f)+Y\rad(B))$ be a projective cover of $Y/({\rm Im}(f)+Y\rad(B))$ and $\pi_2:Y\ra Y/({\rm Im}(f)+Y\rad(B))$ be the natural epimorphism.
Then there is a homomorphism $\gamma_2:Q_B\ra Y_B$ such that $\gamma'_2=\gamma_2\pi_2$.
Then we have the following commutative diagram, and $(P_A,(P\otimes_AM)_B\oplus Q_B,(1_{P\otimes_AM},0)) \stackrel{\gamma}{\ra} (X_A,Y_B,f)_{\Lambda}$
is a projective cover of $(X_A,Y_B,f)_{\Lambda}$, where $\gamma:=(\gamma_1,
\left(\begin{smallmatrix}
(\gamma_1\otimes 1_M)f \\
\gamma_2
\end{smallmatrix}\right)
)$.
$$\xymatrix@C=2cm{
P\otimes_AM
\ar[r]^{(1_{P\otimes_AM},0)}
\ar[d]^{\gamma_1\otimes 1_M}
&(P\otimes_AM)\oplus Q
\ar[d]^{
\left(\begin{smallmatrix}
(\gamma_1\otimes 1_M)f \\
\gamma_2
\end{smallmatrix}\right)
}\\
X\otimes_AM\ar[r]^{f}
&Y
}$$
Let $(\lambda_1,\lambda_2):(W_A,Z_B,\eta)_{\Lambda}\ra (P_A,(P\otimes_AM)_B \oplus Q_B,(1_{P\otimes_AM},0))$ be the kernel of $\gamma$, where $\lambda_1:W_A\ra P_A$ and $\lambda_2:Z_B\ra P\otimes_AM_B\oplus Q_B$. Then we have the following commutative diagram.
$$\xymatrix@R=0.35cm{\\(\star)}
\qquad
\xymatrix@C=2cm{
W\otimes_AM\ar[r]^{\eta}\ar[d]^{\lambda_1\otimes 1_{M}}
&Z\ar[d]^{\lambda_2}\\
P\otimes_AM\ar[r]^{(1_{P\otimes_AM},0)}
&(P\otimes_AM)\oplus Q
}$$
Note that we have an exact sequence
$$0\lra W_A\lraf{\lambda_1} P_A\lraf{\gamma_1} X_A\lra 0 $$
in $\mbox{mod-}A$ and $W_A\simeq \Omega(X_A)$. Since $\gamma_1$ is a projective cover of $X_A$, we have $(W)\lambda_1\subseteq \rad(P_A)$. Then for $\Sigma_jw_j\otimes m_j\in W\otimes_AM$, we have $(\Sigma_jw_j\otimes m_j)(\lambda_1\otimes 1_{M})=\Sigma_j(w_j)\lambda_1\otimes m_j$, where $w_j\in W$ and $m_j\in M$ for all $j$.
Due to $(W)\lambda_1\subseteq \rad(P_A)$, there exist $p_j\in P$ and $r_j\in \rad(A)$ such that $(w_j)\lambda_1=p_jr_j$.
As $\rad(A)M=0$, we have
$$\Sigma_j(w_j)\lambda_1\otimes m_j=\Sigma_jp_jr_j\otimes m_j=\Sigma_jp_j\otimes r_j m_j=0.$$
Thus $\lambda_1\otimes 1_{M}=0$. By the commutative diagram $(\star)$, since $\lambda_2$ is injective, we get $\eta=0$. Then
$(W,Z,\eta)_{\Lambda}\simeq (W,0,0)_{\Lambda}\oplus (0,Z,0)_{\Lambda}$.
 Also $W_A\simeq \Omega(X_A)$.
Then we have the following exact sequence of $\Lambda$-modules:
$$0\lra (\Omega(X_A),0,0)\oplus (0,Z_B,0)\lra (P_A,(P\otimes_AM)_B\oplus Q_B,(1_{P\otimes_AM},0)) \lraf{\gamma} (X_A,Y_B,f)_{\Lambda}\lra 0.$$

Next, we show $Z_B\in \add(B/\rad(B))$. Since $\modcat{B}=\add(B_B\oplus B/\rad(B))$ by Lemma \ref{T(A)}, we have to show that $Z_B$ has no nonzero projective summands. Indeed, assume that $Z_B$ has a nonzero projective summand $Z'_B$. By Lemma \ref{T(A)}, we know that $B$ is a symmetric algebra. In particular, each projective $B$-module is also injective, and $Z'_B$ is injective. Note that we have the following exact sequence of $B$-modules:
$$(*)\quad
\xymatrix@C=1.5cm{
0\ar[r]
&Z_B\ar[r]^(0.32){\lambda_2}
&(P\otimes_AM)_B\oplus Q_B
\ar[r]^(0.65){
\left(\begin{smallmatrix}
(\gamma_1\otimes 1_M)f \\
\gamma_2
\end{smallmatrix}\right)}
&Y_B\ar[r]
&0
}$$
Then $(Z'_B)\lambda_2$ is a direct summand of $(P\otimes_AM)_B\oplus Q_B$. It follows from the definition of $M_B$ that $M_B\simeq \top(B_B)$ and $M\rad(B)=0$. Thus $M_B$ and $(P\otimes_AM)_B$ are semisimple $B$-modules.
By Lemma \ref{T(A)}, $B$ has no projective simple modules, and therefore $(Z'_B)\lambda_2$ is a direct summand of $Q_B$. So, we can write $Q=Q'\oplus (Z')\lambda_2$ for some projective $B$-module $Q'$. It follows from $(*)$ that $((Z')\lambda_2)\gamma_2=0$.
As $\gamma'_2=\gamma_2\pi_2$, we have $((Z')\lambda_2)\gamma'_2=0$.
This is a contradiction since $\gamma'_2:Q_B\ra Y/({\rm Im}(f)+Y\rad(B))$ is a projective cover of $Y/({\rm Im}(f)+Y\rad(B))$. Thus $Z_B$ has no nonzero projective summands, and therefore $Z_B\in \add(B/\rad(B)_B)$.

Let $Z_1:=Z$.
We get the desired exact sequence $(\dagger)$.
In particular, $\Omega((X_A,Y_B,f)_{\Lambda})\simeq (\Omega(X_A),0,0)\oplus (0,Z_1,0)_{\Lambda},\mbox{ and}$
\begin{align*}
\Omega^2((X_A,Y_B,f)_{\Lambda})
&=\Omega(\Omega(X_A,Y_B,f)_{\Lambda})\\
&\simeq \Omega((\Omega(X_A),0,0)_{\Lambda}\oplus (0,Z_1,0)_{\Lambda})
\quad \mbox{(by $(\dagger)$)}\\
&\simeq (\Omega^2(X_A),0,0)\oplus (0,Z'_1,0)\oplus \Omega((0,Z_1,0)_{\Lambda})
\quad \mbox{(by $(\dagger)$)}
\end{align*}
for some $Z'_1\in \add(B/\rad(B))$.
It follows from Lemma \ref{simple-seq} and $Z_1\in \add(B/\rad(B)_B)$ that $\Omega((0,Z_1,0)_{\Lambda})\in add((0,B/\rad(B),0)_{\Lambda})$. Thus $$\Omega^2((X_A,Y_B,f)_{\Lambda})\simeq (\Omega^2(X_A),0,0)\oplus (0,Z_2,0)_{\Lambda}$$
for some $Z_2\in\add(B/\rad(B)_B)$. Continuing this argument, by $(\dagger)$ and  Lemma \ref{simple-seq}, we have $$\Omega^s((X_A,Y_B,f)_{\Lambda})\simeq (\Omega^s(X_A),0,0)_{\Lambda} \oplus (0,Z_s,0)_{\Lambda}$$ with $Z_s\in\add(B/\rad(B)_B)$ for $s\ge 3$.
$\square$

\begin{Theo}\label{tri-del}
$\del(A)\le \del(\Lambda)\le \del(A)+1$.
\end{Theo}
{\rm Proof.}
Let $e:=
\left(\begin{smallmatrix}
1 & 0 \\
0 & 0
\end{smallmatrix}\right)$ and $M:={_AA/\rad(A)_B}$. If $\del(A)=\infty$ or $\del(\Lambda)=\infty$, then the inequalities are trivially true. So, we assume $\del(A)<\infty$ and $\del(\Lambda)<\infty$.

We first show $\del(A)\le \del(\Lambda)$.
Suppose $\del(\Lambda)=d$.
%
%
Then
there exist $(X_A,Y_B,f)\in \modcat{\Lambda}$ and $V\in \pmodcat{\Lambda}$ such that the $d$-th syzygy $\Omega^{d}((A/\rad(A),0,0)_{\Lambda})$ of the semisimple $\Lambda$-module $(A/\rad(A),0,0)$ is a direct summand of $\Omega^{d+1}((X_A,Y_B,f)_{\Lambda})\oplus V$.
By Lemma \ref{tri-del-lem}, we have
$$\Omega^d((A/\rad(A)_A,0,0)_{\Lambda})\simeq (\Omega^d(A/\rad(A)_A),0,0)_{\Lambda} \oplus (0,Z_d,0)_{\Lambda},\mbox{ and}$$
$$\Omega^{d+1}((X_A,Y_B,f)_{\Lambda})\simeq (\Omega^{d+1}(X_A),0,0)_{\Lambda} \oplus (0,Z'_{d+1},0)_{\Lambda}$$
with $Z_d,\, Z'_{d+1}\in\add(B/\rad(B)_B)$.
Then $$(\Omega^d(A/\rad(A)_A),0,0)_{\Lambda}\mbox{ is a direct summand of }(\Omega^{d+1}(X_A),0,0)_{\Lambda} \oplus (0,Z'_{d+1},0)_{\Lambda}\oplus V.$$
Write $V= (P',P'\otimes_AM,1_{P'\otimes_AM})\oplus (0,Q',0)$ for some projective $A$-module $P'$ and some projective $B$-module $Q'$.
It follows from $\Hom_{\Lambda}((N_A,0,0),(0,W_B,0))=0$ for any $N_A$ and $W_B$ that $$(\Omega^d(A/\rad(A)_A),0,0)_{\Lambda}\mbox{ is a direct summand of }(\Omega^{d+1}(X_A),0,0)_{\Lambda} \oplus (P',P'\otimes_AM,1_{P'\otimes_AM})_{\Lambda}.$$
Note that if $P'\neq 0$, then $P'\otimes_AM=P'\otimes_AA/\rad(A)\neq 0$.
Thus $$(\Omega^d(A/\rad(A)_A),0,0)_{\Lambda}\mbox{ is a direct summand of }(\Omega^{d+1}(X_A),0,0)_{\Lambda}.$$ Then
$\Omega^{d}(A/\rad(A)_A)$ is a direct summand of $\Omega^{d+1}(X_A)$ and therefore $\del(A/\rad(A)_A)\le d$. Hence $$\del(A)=\del(A/\rad(A)_A)\le d=\del(\Lambda).$$

Next, we show that $\del(\Lambda)\le \del(A)+1$.
Suppose $\del(A)=n$.
Then there is an $A$-module $X$ such that $\Omega^n(A/\rad(A)_A)$ is a direct summand of $\Omega^{n+1}(X_A)\oplus P$ for some $P\in \pmodcat{A}$.
By Lemma \ref{tri-del-lem}, there exist $Z,Z'\in \add(B/\rad(B))$ such that $$\Omega^{n}((A/\rad(A),0,0)_{\Lambda})\simeq (\Omega^n(A/\rad(A)_A),0,0)\oplus (0,Z_B,0)_{\Lambda},\mbox{ and}$$ $$\Omega^{n+1}((X,0,0)_{\Lambda})\simeq (\Omega^{n+1}(X_A),0,0)\oplus (0,Z'_B,0)_{\Lambda}.$$
Then $$\Omega^{n}((A/\rad(A),0,0)_{\Lambda})\mbox{ is a direct summand of }\Omega^{n+1}((X,0,0)_{\Lambda})\oplus (P_A,0,0)_{\Lambda}\oplus (0,Z_B,0)_{\Lambda}.$$
It follows from $\Omega((P,0,0)_{\Lambda})\simeq (0,P\otimes_AM,0)_{\Lambda}$ that $$\Omega^{n+1}((A/\rad(A),0,0)_{\Lambda})\mbox{ is a direct summand of }\Omega^{n+2}((X,0,0)_{\Lambda})\oplus (0,P\otimes_AM,0)_{\Lambda}\oplus \Omega((0,Z_B,0)_{\Lambda}).$$
By Lemma \ref{simple-seq}, $\add(\Omega^i((0,B/\rad(B),0)_{\Lambda}))=\add((0,B/\rad(B),0)_{\Lambda})$ for $i\ge 0$.
Due to $P\otimes_AM, Z\in \add(B/\rad(B)_B)$, there exists $Z_0\in B/\rad(B)_B$ such that $(0,P\otimes_AM,0)_{\Lambda}\oplus \Omega((0,Z_B,0)_{\Lambda})$ is a direct summand of $\Omega^{n+2}((0,Z_0,0)_{\Lambda})$.
Thus $$\Omega^{n+1}((A/\rad(A),0,0)_{\Lambda})\mbox{ is a direct summand of }\Omega^{n+2}((X,0,0)_{\Lambda}\oplus (0,Z_0,0)_{\Lambda}),$$ and $\del((A/\rad(A),0,0)_{\Lambda})\le n+1$.
It follows from Lemma \ref{simple-seq} that $\del((0,B/\rad(B),0)_{\Lambda})=0$.
Hence $\del(\Lambda)=\max\{
\del((A/\rad(A),0,0)_{\Lambda}),
\del((0,B/\rad(B),0)_{\Lambda})
\}\le n+1=\del(A)+1$.
$\square$

Now, we display examples to illustrate the constructions in Constructions \ref{consruction1} and \ref{consruction2}, and to show that the two cases $\del(\Lambda)=\del(A)$ and $\del(\Lambda)=\del(A)+1$ in Theorem \ref{tri-del} can really happen.
\begin{Bsp}{\rm
(1) Let $A$ be the algebra (over a field) given by the following quiver with a relation:
$$\xymatrix{
1\ar@<0.5ex>[r]^{\alpha} & 2,\ar@<0.5ex>[l]^{\beta}&\beta\alpha=0
}$$
(see \cite[Example 4.9]{Xi2021} for the general case). It follows from \cite[Theorem 1.1]{Xi2022} and \cite[Example 1.18]{g22} that $A$ is Gorenstein and $\del(A)=2$.
The Auslander-Reiten quiver of $A$ can be seen as follows:
$$\xymatrix@R=0.8cm@C=0.4cm{
\ar@{.}[d]&&&&\ar@{.}[d]\\
S_A(1)\ar[dr]\ar@{.}[d]
&&S_A(2)\ar@{-->}[ll]\ar[dr]
&&S_A(1)\ar@{-->}[ll]\ar@{.}[d]\\
\ar@{.}[d]
&
P_A(2)
\ar[ur]\ar[dr]
&&
I_A(2)
\ar@{-->}[ll]\ar[ur]
&\ar@{.}[d]\\
&&P_A(1)\ar[ur]
&&
}$$
where $S_A(i)$, $P_A(i)$ and $I_A(i)$ are the simple, projective and injective $A$-modules corresponding to each vertex $i$, respectively.
Here the vertical dotted lines are identified.

Let $\widetilde{A}$ and $\Lambda$ be defined as in Construction \ref{consruction1} and \ref{consruction2}, respectively. By calculation, $\widetilde{A}$ and $\Lambda$ can be described by the quivers with relations, respectively.
$$\begin{array}{ccc}
\widetilde{A}:\;\xymatrix{
\tilde{1}\ar@(ur,ul)_{\sigma_1}\ar[d]^{\tau_1}
&\tilde{2}\ar@(ur,ul)_{\sigma_2}\ar[d]^{\tau_2}\\
1 \ar@<0.5ex>[r]^{\alpha}
&2\ar@<0.5ex>[l]^{\beta}}\\
\sigma_1^2=\sigma_2^2=\sigma_1\tau_1=\sigma_2\tau_2=\tau_1\alpha=\tau_2\beta=\beta\alpha=0,
\end{array}
\quad
\begin{array}{ccc}
\Lambda:\;\xymatrix{
1 \ar@<0.5ex>[r]^{\alpha}\ar[d]^{\epsilon_1}
&2\ar@<0.5ex>[l]^{\beta}\ar[d]^{\epsilon_2}\\
1'\ar@(dr,dl)^{\eta_1}
&2'\ar@(dr,dl)^{\eta_2}}\\
\eta_1^2=\eta_2^2=\epsilon_1\eta_1=\epsilon_2\eta_2=\beta\epsilon_1=\alpha\epsilon_2=\beta\alpha=0
\end{array}$$
The Loewy structures of indecomposable projective $\widetilde{A}$-modules and $\Lambda$-modules are displayed, respectively.
$$
\xymatrix@R=.23cm@C=.01cm{
e_1\widetilde{A}\\
1\ar@{-}[d]\\
2\ar@{-}[d]\\
1}\quad
\xymatrix@R=.23cm@C=.01cm{
e_2\widetilde{A}\\
2\ar@{-}[d]\\
1}\quad
\xymatrix@R=.23cm@C=.01cm{
&e_{\tilde{1}}\widetilde{A}&\\
&\tilde{1}\ar@{-}[dl]\ar@{-}[dr]&\\
\tilde{1}&&1}\quad
\xymatrix@R=.23cm@C=.01cm{
&e_{\tilde{2}}\widetilde{A}&\\
&\tilde{2}\ar@{-}[dl]\ar@{-}[dr]&\\
\tilde{2}&&2}
\qquad
\xymatrix@R=.23cm@C=.01cm{
&e_1\Lambda&\\
&1\ar@{-}[dl]\ar@{-}[dr]&\\
1'&&2\ar@{-}[d]\\
&&1}\quad
\xymatrix@R=.23cm@C=.01cm{
&e_2\Lambda&\\
&2\ar@{-}[dl]\ar@{-}[dr]&\\
2'&&1}\quad
\xymatrix@R=.23cm@C=.01cm{
e_{1'}\Lambda\\
1'\ar@{-}[d]\\
1'}\quad
\xymatrix@R=.23cm@C=.01cm{
e_{2'}\Lambda\\
2'\ar@{-}[d]\\
2'}
$$
Clearly, each simple $\widetilde{A}$-module is a submodule of some projective $\widetilde{A}$-module. Thus $\del(\widetilde{A})=0$.

Next, we show that $\del(\Lambda)=3$. Indeed, let $S_{\Lambda}(i)$ be the simple $\Lambda$-module corresponding to each vertex $i$. Then we have the following exact sequences of $\Lambda$-modules.
$$0\lra S_{\Lambda}(2')\oplus S_{\Lambda}(1) \lra e_{2}\Lambda \lra S_{\Lambda}(2)\lra 0$$
$$(*)\quad 0\lra S_{\Lambda}(1') \lra e_{1'}\Lambda \lra S_{\Lambda}(1')\lra 0$$
$$(**)\quad 0\lra S_{\Lambda}(2') \lra e_{2'}\Lambda \lra S_{\Lambda}(2')\lra 0$$
Thus $S_{\Lambda}(1)\oplus S_{\Lambda}(1')\oplus S_{\Lambda}(2')\in \add(S_{\Lambda}(2)\oplus S_{\Lambda}(1')\oplus S_{\Lambda}(2'))$ and $\del(S_{\Lambda}(1)\oplus S_{\Lambda}(1')\oplus S_{\Lambda}(2'))=0$. Clearly, $\del(S_{\Lambda}(2))\ge 1$. Hence $\del(\Lambda)=\del(S_{\Lambda}(2))\ge 1$.

Note that we have the following exact sequences of $\Lambda$-modules.
$$0\lra S_{\Lambda}(1')\oplus P_2(A) \lra e_{1}\Lambda \lra S_{\Lambda}(1)\lra 0$$
$$0\lra S_{\Lambda}(2') \lra e_{2}\Lambda \lra P_2(A)\lra 0$$
Here we view $P_2(A)$ as a $\Lambda$-module.
Then $$\Omega(S_{\Lambda}(2))\simeq S_{\Lambda}(2')\oplus S_{\Lambda}(1),$$ $$\Omega^2(S_{\Lambda}(2))\simeq S_{\Lambda}(2')\oplus S_{\Lambda}(1')\oplus P_2(A),\mbox{ and}$$ $$\Omega^3(S_{\Lambda}(2))\simeq S_{\Lambda}(2')\oplus S_{\Lambda}(1')\oplus S_{\Lambda}(2').$$
By $(*)$ and $(**)$, we have $\Omega^i(S_{\Lambda}(1'))\simeq S_{\Lambda}(1')$ and $\Omega^i(S_{\Lambda}(2'))\simeq S_{\Lambda}(2')$ for $i\ge 0$. Thus $$\Omega^3(S_{\Lambda}(2))\simeq \Omega^4(S_{\Lambda}(2')\oplus S_{\Lambda}(1')\oplus S_{\Lambda}(2'))$$ and $\del(S_{\Lambda}(2))\le 3$.

By Lemma \ref{tri-del-lem}, $\Omega^s(\modcat{\Lambda})=\add(\{(X,0,0)\mid X\in \Omega^s(\modcat{A}) \}\oplus S_{\Lambda}(1')\oplus S_{\Lambda}(2'))$, where $\Omega^s(\modcat{\Lambda})$ stands for the full subcategory of $\modcat{\Lambda}$ consisting of direct summands of $s$-syzygies of $\Lambda$-modules for $s\ge 1$.
By calculation, we get $\Omega^2(\modcat{A})=\add(P_A(2))$ and $\Omega^3(\modcat{A})=0$.
Then 
$$\Omega^2(\modcat{\Lambda})=\add((P_A(2),0,0)_{\Lambda}\oplus S_{\Lambda}(1')\oplus S_{\Lambda}(2'))\mbox{ and }
\Omega^3(\modcat{\Lambda})=\add(S_{\Lambda}(1')\oplus S_{\Lambda}(2')).$$
Thus $\Omega(S_{\Lambda}(2))\not\in \Omega^2(\modcat{\Lambda})\oplus \add(\Lambda_{\Lambda})\mbox{ and }\Omega^2(S_{\Lambda}(2))\not\in \Omega^3(\modcat{\Lambda})\oplus \add(\Lambda_{\Lambda}).$
So $\del(S_{\Lambda}(2))\neq 1,2$.
It follows $\del(S_{\Lambda}(2))\ge 1$ that $\del(S_{\Lambda}(2))\ge 3$, and therefore $\del(\Lambda)=\del(S_{\Lambda}(2))= 3=\del(A)+1$.

(2) Let $A$ be the algebra (over a field) given by the following quiver with a relation:
$$\xymatrix{
1\ar@(ru,rd)^{\alpha}&\qquad\alpha^2=0.
}$$
Since the unique simple module $S_A(1)$ of $A$ is $1$-syzygy of $S_A(1)$, we obtain $\del(A)=\del(S_A(1))=0$.

Let $\widetilde{A}$ and $\Lambda$ be defined as in Construction \ref{consruction1} and \ref{consruction2}, respectively. By calculation, $\widetilde{A}$ and $\Lambda$ can be described by the quivers with relations, respectively.
$$\begin{array}{ccc}
\widetilde{A}:\;\xymatrix{
\tilde{1}\ar@(ru,rd)^{\sigma}\ar[d]^{\tau}\\
1\ar@(ru,rd)^{\alpha}
}\\
\\
\sigma^2=\sigma\tau=\tau\alpha=\alpha^2=0
\end{array}
\qquad
\begin{array}{ccc}
\Lambda:\;\xymatrix{
1\ar@(ru,rd)^{\alpha}\ar[d]^{\epsilon}\\
1'\ar@(ru,rd)^{\eta}
}\\
\\
\eta^2=\epsilon\eta=\alpha\epsilon=\alpha^2=0
\end{array}
$$
The Loewy structures of indecomposable projective $\Lambda$-modules are displayed as follows:
$$
\xymatrix@R=.23cm@C=.01cm{
&e_1\Lambda&\\
&1\ar@{-}[dl]\ar@{-}[dr]&\\
1&&1'}\quad
\xymatrix@R=.23cm@C=.01cm{
e_2\Lambda\\
1'\ar@{-}[d]\\
1'}
$$
Clearly, $\widetilde{A}\simeq \Lambda$ as algebras, and $\del(\widetilde{A})=\del(\Lambda)=0$. Thus $\del(\Lambda)=\del(A)$.
}\end{Bsp}

As an application of Theorem \ref{main-thm}, we have the following consequence. 
\begin{Prop}\label{left-right}
The delooping level of any Artin algebra is finite if and only if, for all Artin algebras, the delooping level of an algebra being finite implies that the delooping level of its opposite algebra is also finite.
\end{Prop}
{\rm Proof.} The forward direction holds trivially, so we focus on the reverse direction.

Let $A$ be an Artin algebra and $\Lambda$ be defined as in Construction \ref{consruction2}. Then $\del(\Lambda^{\opp})=0$ by Lemma \ref{delopp=0}. Hence, by assumption,
$$\del(\Lambda)=\del((\Lambda^{\opp})^{\opp}) < \infty.$$
It follows from Lemma \ref{tri-del} that $\del(A)\le \del(\Lambda) <\infty$.
$\square$

Finally, we present an example of a finite-dimensional algebra with an infinite delooping level, while its opposite algebra has a delooping level of zero, based on an example by Kershaw and Rickard \cite{kr23}.

\begin{Bsp}{\rm
Let $k$ be a field,  and $q\in k$ an element with infinite multiplicative order. Let $C$ be the algebra
$$k<x,y,z>/(x^2,y^2,z^2,zy,yx+qxy,zx-xz,yz-xz),$$
and  $M$ the three-dimensional $C$-module with basis 
$v,v',v''$, where $vx=qv',vy=v',vz=v''$ with $v'$ and $v''$ annihilated by $x,y$ and $z$. By  \cite[Corollary 4.2]{kr23}, the one-point extension algebra 
$$C[M]:=\begin{pmatrix}
k&M\\
0&C
\end{pmatrix}$$ 
has infinite delooping level. 
Denote by $\Lambda$ the related algebra of $A:=C[M]$ as Construction \ref{consruction2}. By Theorem \ref{main-thm}, we have $\del(\Lambda)=\infty$ and $\del(\Lambda^{\opp})=0$.
}\end{Bsp}

\bigskip
\noindent{\bf Acknowledgements.}
This work was supported by the National Natural Science Foundation of China (Grants 12501052 and 12401038).

\noindent{\bf Declaration of interests.} The authors have no conflicts of interest to disclose.

\noindent{\bf Data availability.} No new data were created or analyzed in this study.

\end{document}